# Efficient three-stage $t$-tests

## Jay Bartroff[1],*

*Stanford University*

**Abstract:** Three-stage $t$-tests of separated one-sided hypotheses are derived, extending Lorden's optimal three-stage tests for the one-dimensional exponential family by using Lai and Zhang's generalization of Schwarz's optimal fully-sequential tests to the multiparameter exponential family. The resulting three-stage $t$-tests are shown to be asymptotically optimal, achieving the same average sample size as optimal fully-sequential tests.

## 1. Introduction

Lorden [5] derived asymptotically optimal three-stage tests of separated one-sided hypotheses about the parameter of a one-dimensional exponential family. Lorden's tests use the initial stage to estimate the unknown parameter and then use this to choose the size of the second stage close to the average sample size of Schwarz's [6] optimal sequential tests. Lorden showed that the resulting three-stage tests are as asymptotically efficient as optimal fully-sequential tests, asymptotically achieving Hoeffding's [1] lower bound on the average sample size of a sequential test with given error probabilities. Lorden [5] showed that, conversely, three stages were also necessary for asymptotic efficiency, except in a few degenerate cases. Moreover, Lorden [5] showed that any efficient three-stage test must mimic Schwarz's sequential test in this way.

Lai [2] modified Schwarz's [6] tests of separated one-sided hypotheses to allow the distance between the hypotheses to approach zero, and hence derived asymptotically optimal sequential tests of hypotheses with or without indifference regions. Lai and Zhang [4] extended Lai's [2] tests to general hypotheses in the multiparameter exponential family setting by developing and applying certain results on boundary crossing probabilities of multiparameter generalized likelihood ratio statistics.

In this paper we present three-stage $t$-tests of one-sided, separated hypotheses about the mean of i.i.d. normal data whose variance is unknown. On one hand, this is an extension of the tests of Lorden [5] to the two-parameter setting. On the other hand, this is a first step in adapting the multiparameter sequential tests of Lai and Zhang [4] to the multistage setting. Theorem 2.2 shows that the expected sample size of our three-stage $t$-test asymptotically achieves the Hoeffding [1] lower bound as the error probabilities approach zero. Converse results and extensions are discussed in Section 3.

## 2. Three-stage $t$-tests

Let $X_1, X_2, \ldots$ be i.i.d. $N(\mu, \sigma^2)$, where both $\mu$ and $\sigma^2$ are unknown. We consider testing

---

[1]Department of Statistics, Sequoia Hall, Stanford University, Stanford, CA 94305, e-mail: `bartroff@stat.stanford.edu`

*This work is supported by NSF Grant DMS-0403105

*AMS 2000 subject classifications:* primary 62F05; secondary 62L10.

*Keywords and phrases:* multistage hypothesis test, $t$-test, asymptotic efficiency.

$$H_0 : \mu \leq \mu_0 \quad \text{vs.} \quad H_1 : \mu \geq \mu_1 > \mu_0.$$

It is assumed that

$$0 < \underline{\sigma}^2 \leq \sigma^2 \leq \overline{\sigma}^2.$$

In practice, such bounds may be indicated by prior experience with the type of data or they may be implied by practical considerations like desired maximum and minimum total sample sizes. For example, if the $X_i$ represent the difference between paired treatment and control responses in a clinical trial, previous experience with similar trials may suggest values of $\underline{\sigma}^2, \overline{\sigma}^2$. Alternatively, it may be more natural for practitioners to work in terms of minimum and maximum sample size, which $\underline{\sigma}^2$ and $\overline{\sigma}^2$ imply, as will be seen below.

Let $L_n(u, v^2)$ denote the log-likelihood function

$$L_n(u, v^2) = -(n/2) \log v^2 - (2v^2)^{-1} \sum_{i=1}^n (X_i - u)^2.$$

The log of the generalized likelihood ratio (GLR) is, for $i = 0, 1$,

$$\Lambda_{i,n} = \sup_{u,v^2} L_n(u, v^2) - \sup_{v^2} L_n(\mu_i, v^2)$$

$$(2.1) \qquad\qquad = (n/2) \log \left[ 1 + \left( \frac{\overline{X}_n - \mu_i}{\widehat{\sigma}_n} \right)^2 \right]$$

$$(2.2) \qquad\qquad = n I_i(\overline{X}_n, \widehat{\sigma}_n^2),$$

where $\overline{X}_n$ is the sample mean,

$$\widehat{\sigma}_n^2 = \frac{1}{n} \sum_{i=1}^n (X_i - \overline{X}_n)^2, \quad \text{and}$$

$$I_i(u, v^2) = (1/2) \log[1 + (u - \mu_i)^2 / v^2]$$

is the Kullback–Leibler information number. The term being squared inside the log in (2.1) is of course a multiple of the usual $t$ statistic, but it will be more natural in what follows to work with $\Lambda_{i,n}$ rather than the $t$ statistic. Given $A_0, A_1 > 0$, a GLR test of $H_0$ vs. $H_1$ rejects $H_0$ if

$$(2.3) \qquad\qquad \overline{X}_n > \mu_0 \quad \text{and} \quad \Lambda_{0,n} \geq A_0,$$

and rejects $H_1$ if

$$(2.4) \qquad\qquad \overline{X}_n < \mu_1 \quad \text{and} \quad \Lambda_{1,n} \geq A_1.$$

The boundaries (2.3) and (2.4) define "stopping surfaces" in $(n, \sum_1^n X_i, \sum_1^n (X_i - \overline{X}_n)^2)$-space. The surface (2.3) intersects the line

$$(2.5) \qquad\qquad \sum_1^n X_i = un, \quad \sum_1^n (X_i - \overline{X}_n)^2 = v^2 n$$

at a point whose $n$-coordinate is

$$n_0(u, v^2) = \frac{A_0}{I_0(u, v^2)},$$

where we let $n_0(\mu_0, v^2) = \infty$ since $I_0(\mu_0, v^2) = 0$. Define $n_1(u, v^2) = A_1/I_1(u, v^2)$ similarly with respect to the surface (2.4), and let

$$n(u, v^2) = \min_{i=0,1} n_i(u, v^2).$$

Let $\mu_2 = \mu_2(A_0, A_1) \in (\mu_0, \mu_1)$ be the unique solution of

$$\frac{A_0}{A_1} = \frac{I_0(\mu_2, \overline{\sigma}^2)}{I_1(\mu_2, \overline{\sigma}^2)},$$

which exists because the right hand side increases from $0$ to $+\infty$ as $\mu_2$ ranges from $\mu_0$ to $\mu_1$, and let $\overline{n} = n(\mu_2, \overline{\sigma}^2)$. Since $n_i(u, v^2) \leq n_i(u, \overline{\sigma}^2)$ for all $u$ and $v^2 < \overline{\sigma}^2$,

$$(2.6) \qquad \overline{n} = n(\mu_2, \overline{\sigma}^2) = \sup_u n(u, \overline{\sigma}^2) = \sup_{u, v^2 < \overline{\sigma}^2} n(u, v^2).$$

Hence $\overline{n}$ is roughly the "worst-case" expected total sample size required to cross one of the boundaries (2.3) or (2.4) when the variance is no greater than $\overline{\sigma}^2$.

Asymptotic optimality results will be proved for bounded values of $\mu \in [\underline{\mu}, \overline{\mu}]$ where

$$\underline{\mu} < \mu_0 < \mu_1 < \overline{\mu}.$$

As with the bounds on $\sigma^2$, in practice these values may be suggested by prior knowledge or by the constraints on maximum and minimum sample size that they imply. Let

$$J = [\underline{\mu}, \overline{\mu}] \times [\underline{\sigma}^2, \overline{\sigma}^2],$$
$$J_\varepsilon = (\underline{\mu} - \varepsilon, \overline{\mu} + \varepsilon) \times (\underline{\sigma}^2 - \varepsilon, \overline{\sigma}^2 + \varepsilon)$$

for $\varepsilon > 0$.

Before defining the stopping rule of our three stage $t$-test we state an auxiliary lemma that gives bounds on how close the estimated total sample size after sampling $k$, $n(\overline{X}_k, \widehat{\sigma}_k^2)$, is to the "correct" total sample size $n(\mu, \sigma^2)$. This is an extension of Lemma 1 of Lorden [5] to the two parameter setting.

**Lemma 2.1.** *Let* $0 < \varepsilon < \underline{\sigma}^2$. *There is a positive constant* $B$ *such that if*

$$\rho_n = 1 + B\sqrt{n^{-1} \log n},$$

*then*

$$(2.7) \qquad \begin{aligned} P_{\mu, \sigma^2} &\left( (\overline{X}_k, \widehat{\sigma}_k^2) \in J_\varepsilon, \rho_k^{-1} < \frac{n(\overline{X}_k, \widehat{\sigma}_k^2)}{n(\mu, \sigma^2)} < \rho_k, \text{for all } k \geq k_0 \right) \\ &\geq 1 - O(1/k_0) \end{aligned}$$

*as* $k_0 \to \infty$, *uniformly for* $(\mu, \sigma^2) \in J$.

*Proof.* Using large deviations probabilities (e.g., Lai & Zhang [4], Lemma 1) it can be shown that there is a $\beta > 0$ such that

$$(2.8) \qquad P_{\mu, \sigma^2}\left( |\overline{X}_k - \mu| \geq \eta \quad \text{some } k \geq k_0 \right) \leq \exp(-\beta\eta^2 k_0^2)$$
$$(2.9) \qquad P_{\mu, \sigma^2}\left( |\widehat{\sigma}_k^2 - \sigma^2| \geq \eta \quad \text{some } k \geq k_0 \right) \leq \exp(-\beta\eta^2 k_0^2)$$

as $k_0 \to \infty$, uniformly in $(\mu, \sigma^2) \in J_\varepsilon$. By a straightforward Taylor series argument it can be shown that there is $D < \infty$ such that

$$(2.10) \qquad \left| \frac{n(u, v^2)}{n(\mu, \sigma^2)} - 1 \right| \leq D(|u - \mu| + |v^2 - \sigma^2|)$$

for all $(u, v^2) \in J_\varepsilon$, $(\mu, \sigma^2) \in J$.

Let $V$ be the event in (2.7). Denote $P_{\mu, \sigma^2}$ simply by $P$. First write

$$
\begin{aligned}
(2.11) \quad P(V^c) \leq{}& P\left( (\overline{X}_k, \widehat{\sigma}_k^2) \in J_\varepsilon, \left| \frac{n(\overline{X}_k, \widehat{\sigma}_k^2)}{n(\mu, \sigma^2)} - 1 \right| > \frac{\rho_k - 1}{2}, \text{some } k \geq k_0 \right) \\
&+ P((\overline{X}_k, \widehat{\sigma}_k^2) \notin J_\varepsilon, \text{some } k \geq k_0).
\end{aligned}
$$

The latter is no greater than

$$P(|\overline{X}_k - \mu| \geq \varepsilon, \text{some } k \geq k_0) + P(|\widehat{\sigma}_k^2 - \sigma^2| \geq \varepsilon, \text{some } k \geq k_0) \leq 2\exp(-2\beta' k_0^2)$$

for some $\beta' > 0$ by (2.8) and (2.9). Since this is obviously $o(1/k_0)$, the proof of the lemma will be complete once we show that the former term in (2.11) is $O(1/k_0)$. Using (2.10), this term is no greater than the sum over all $k \geq k_0$ of

$$
\begin{aligned}
&P\left( |\overline{X}_k - \mu| + |\widehat{\sigma}_k^2 - \sigma^2| > \frac{\rho_k - 1}{2D} \right) \\
&\qquad \leq P\left( |\overline{X}_k - \mu| > \frac{\rho_k - 1}{4D} \right) + P\left( |\widehat{\sigma}_k^2 - \sigma^2| > \frac{\rho_k - 1}{4D} \right).
\end{aligned}
$$

By (2.8) and (2.9), for some $\beta'' > 0$ this is no greater than

$$\sum_{k \geq k_0} 2\exp[-\beta''(\rho_k - 1)^2 k^2] = \sum_{k \geq k_0} 2\exp(-\beta'' B \log k) = \sum_{k \geq k_0} 2/k^2 = O(1/k_0),$$

by choosing $B = 2/\beta''$.                                                                                    □

We can now define the stopping times $N_1, N_2, N_3$ of the three stage $t$-test. Let $m = \lceil \overline{n} \rceil$, $A_0, A_1 > 0$, and choose $0 < C \leq 1$ and $0 < \varepsilon < \underline{\sigma}^2$. Let

$$
\begin{aligned}
N_1 &= \lfloor C[n(\overline{\mu}, \underline{\sigma}^2) \wedge n(\underline{\mu}, \underline{\sigma}^2)] \rfloor \\
N_2 &= m \wedge (N_1 \vee \lceil \rho_{N_1}^2 n(\overline{X}_{N_1}, \widehat{\sigma}_{N_1}^2) \rceil) \\
N_3 &= m.
\end{aligned}
$$

The test stops at the end of stage $i \leq 3$ and rejects $H_0$ if

$$(2.12) \qquad \underline{\sigma}^2 - \varepsilon < \widehat{\sigma}_{N_i}^2 < \overline{\sigma}^2 + \varepsilon, \quad \mu_0 < \overline{X}_{N_i} < \overline{\mu} + \varepsilon, \quad \text{and} \quad \Lambda_{0, N_i} \geq A_0.$$

The test stops at the end of stage $i \leq 3$ and rejects $H_1$ if

$$(2.13) \qquad \underline{\sigma}^2 - \varepsilon < \widehat{\sigma}_{N_i}^2 < \overline{\sigma}^2 + \varepsilon, \quad \underline{\mu} - \varepsilon < \overline{X}_{N_i} < \mu_1, \quad \text{and} \quad \Lambda_{1, N_i} \geq A_1.$$

Note that when the total sample size is $m$ and $\widehat{\sigma}_m^2 \leq \overline{\sigma}^2$, one of (2.12), (2.13) is guaranteed to happen by our choice of $\overline{n}$ in (2.6). In the event that the total sample size is $m$ but neither (2.12) nor (2.13) happens, we define the test to reject $H_0$ if and only if $\Lambda_{0, m} > \Lambda_{1, m}$.

**Theorem 2.2.** *Let $N = N(A_0, A_1)$ be the total sample size of the test defined above, and let*

$$(2.14) \qquad \alpha_i = P_{\mu_i, \overline{\sigma}^2}(Reject\ H_i) = \sup_{H_i} P_{\mu, \sigma^2}(Reject\ H_i) \quad (i = 0, 1).$$

*Let $\mathcal{C}(\alpha_0, \alpha_1)$ be the class of all (possibly sequential) tests of $H_0$ vs. $H_1$ with error probabilities (2.14) no greater than $\alpha_0, \alpha_1$. As $A_0, A_1 \to \infty$ such that $A_0 \sim A_1$,*

$$(2.15) \qquad \begin{aligned} E_{\mu, \sigma^2} N &\le n(\mu, \sigma^2) + O(\sqrt{A_i \log A_i}) \\ &= \inf_{N' \in \mathcal{C}(\alpha_0, \alpha_1)} E_{\mu, \sigma^2} N' + O(\sqrt{A_i \log A_i}) \end{aligned}$$

*uniformly for $(\mu, \sigma^2) \in J$, for either $i = 0, 1$.*

*Proof.* To bound $E_{\mu, \sigma^2} N$ from above we follow the three-stage argument of Lorden [5] with Lemma 2.1. Let $V$ be the event in (2.7) with $k_0 = N_1$. We will show that

$$(2.16) \qquad N \le \lceil \rho_{N_1}^2 n(\overline{X}_{N_1}, \hat{\sigma}_{N_1}^2) \rceil \quad \text{on } V.$$

It is true that $N_1 \le n(\mu, \sigma^2)$ since in fact

$$(2.17) \qquad \begin{aligned} n(\mu, \sigma^2) &\ge \inf_{(u, v^2) \in J} n(\mu, \sigma^2) \\ &= \inf_u n(u, \underline{\sigma}^2) \\ &= n(\underline{\mu}, \underline{\sigma}^2) \wedge n(\overline{\mu}, \underline{\sigma}^2) \\ &\ge C^{-1} N_1 \ge N_1, \end{aligned}$$

where (2.17) holds since $n_0(\cdot, \underline{\sigma}^2)$ is non-decreasing and $n_1(\cdot, \underline{\sigma}^2)$ is non-increasing. If the test stops after the first stage,

$$N = N_1 \le n(\mu, \sigma^2) \le \rho_{N_1} n(\overline{X}_{N_1}, \hat{\sigma}_{N_1}^2) \le \lceil \rho_{N_1}^2 n(\overline{X}_{N_1}, \hat{\sigma}_{N_1}^2) \rceil$$

on $V$. If $N = N_2 = m$, then by definition of $N_2$,

$$N = m \le \lceil \rho_{N_1}^2 n(\overline{X}_{N_1}, \hat{\sigma}_{N_1}^2) \rceil$$

on $V$. Otherwise, $N_2 < m$ and so, on $V$,

$$\begin{aligned} N_2 &= \lceil \rho_{N_1}^2 n(\overline{X}_{N_1}, \hat{\sigma}_{N_1}^2) \rceil \\ &\ge \rho_{N_1}^2 n(\overline{X}_{N_1}, \hat{\sigma}_{N_1}^2) \\ &\ge \rho_{N_1} n(\mu, \sigma^2) \\ &\ge (\rho_{N_1} / \rho_{N_2}) n(\overline{X}_{N_2}, \hat{\sigma}_{N_2}^2) \ge n(\overline{X}_{N_2}, \hat{\sigma}_{N_2}^2). \end{aligned}$$

This, along with $(\overline{X}_{N_2}, \hat{\sigma}_{N_2}^2) \in J_\varepsilon$, implies that two stages suffice for stopping in this case. Hence

$$N = N_2 \le \lceil \rho_{N_1}^2 n(\overline{X}_{N_1}, \hat{\sigma}_{N_1}^2) \rceil,$$

establishing (2.16).

Conditioning on $V$ and using Lemma 2.1 gives

$$(2.18) \qquad \begin{aligned} E_{\mu, \sigma^2} N &\le \lceil \rho_{N_1}^3 n(\mu, \sigma^2) \rceil + O(m/N_1) \\ &\le n(\mu, \sigma^2) + O(m(\rho_{N_1} - 1)) + O(1) \\ &= n(\mu, \sigma^2) + O(\sqrt{A_i \log A_i}), \end{aligned}$$

since

$$m(\rho_{N_1} - 1) = mB\sqrt{N_1^{-1}\log N_1}$$
$$= O\left(A_i\sqrt{(\log A_i)/A_i}\right)$$
$$= O\left(\sqrt{A_i\log A_i}\right).$$

Consider the fully-sequential test with rejection rule (2.12) and (2.13) with $N_i$ replaced by $n$, $N_1 \le n \le m$. Lai [3] established that the error probabilities of this test are $O(e^{-A_i}A_i^2)$. Since the error probabilities of this fully sequential test obviously bound the error probabilities of our three-stage test, we have

$$\alpha_i = O(e^{-A_i}A_i^2).$$

Using this estimate with Hoeffding's [1] lower bound in the two-parameter exponential family setting gives

$$\inf_{N' \in \mathcal{C}(\alpha_0, \alpha_1)} E_{\mu,\sigma^2} N' \ge \frac{\log(\alpha_0 + \alpha_1)^{-1}}{\max_j I_j(\mu,\sigma^2)} - O\left(\sqrt{\log(\alpha_0 + \alpha_1)^{-1}}\right)$$
$$\ge \min_j \left\{\frac{\log \alpha_j^{-1}}{I_j(\mu,\sigma^2)}\right\} - O\left(\sqrt{\log(\alpha_0 + \alpha_1)^{-1}}\right)$$
$$= \min_j \left\{\frac{A_j}{I_j(\mu,\sigma^2)}\right\} - O\left(\sqrt{A_i}\right) \quad (i = 0, 1)$$
$$= n(\mu,\sigma^2) - O\left(\sqrt{A_i}\right).$$

Combining this with our bound for $E_{\mu,\sigma^2} N$ in (2.18) gives the desired result.    □

## 3. Converse results and extensions

It is natural to suspect that the converse results obtained by Lorden [5] hold in this two-parameter setting. Namely, that essentially no two-stage test can achieve the asymptotic efficiency (2.15) and that any efficient three-stage test must have a total sample size close to $n(\mu,\sigma^2)$ after it's second stage. Indeed, it is easily shown that Theorem 2 and Corollary 2 of Lorden [5] extend immediately to this setting. However, the $t$-test may be deceptively simple since the hypotheses in question are one-dimensional. It is not yet clear whether these properties of optimal tests will carry over to the general multiparameter setting, especially when the dimensions of the hypotheses exceed one or there is no type of "separation" between the hypotheses, as we had here. A first step toward answering these questions would be to adapt the asymptotically optimal fully-sequential tests of Lai & Zhang [4] to the general multistage setting, which remains to be done.

## Acknowledgments

The author thanks Tze Lai for helpful discussions on this subject.

## References

[1] HOEFFDING, W. (1960). Lower bounds for the expected sample size and the average risk of a sequential procedure. *Ann. Math. Statist.* **31** 352–368. MR0120750

[2] LAI, T. L. (1988). Nearly optimal sequential tests of composite hypotheses. *Ann. Statist.* **16** 856–886. MR0947582

[3] LAI, T. L. (1988). Boundary crossing problems for sample means. *Ann. Probab.* **16** 375–396. MR0920279

[4] LAI, T. L. AND ZHANG, L. (1994). A modification of Schwarz's sequential likelihood ratio tests in multivariate sequential analysis. *Seq. Analysis* **13** 79–96. MR1278468

[5] LORDEN, G. (1983). Asymptotic efficiency of three-stage hypothesis tests. *Ann. Statist.* **11** 129–140. MR0684871

[6] SCHWARZ, G. (1962). Asymptotic shapes of Bayes sequential testing regions. *Ann. Math. Statist.* **33** 224–236. MR0137226